\documentclass{qjmam}
\usepackage{amsmath,amsfonts,amssymb,amsbsy}
\begin{document} 

\title[{Asymptotic analysis of the substrate effect}]{Asymptotic analysis of the substrate effect for an arbitrary indenter}
\author[I.~I.~Argatov and F.~J.~Sabina]{I. I. ARGATOV and F.~J.~SABINA\footnote{Corresponding author}}
\address{Institute of Mathematics and Physics, Aberystwyth University, \\ Ceredigion SY23 3BZ, Wales, UK}
\extraaddress{Instituto de Investigaciones en Matem\'aticas
Aplicadas y en Sistemas, Universidad Nacional Aut\'onoma de
M\'exico, Apartado Postal 20-726, Delegaci\'on Alvaro Obregon
01000, M\'exico, D.~F., M\'exico}
\label{firstpage}

\received{\recd %25 March 1997.
\revd
% 24 September 1997
}
\maketitle
\eqnobysec
\begin{abstract}
A quasistatic unilateral frictionless contact problem for a rigid axisymmetric indenter pressed into a homogeneous, linearly elastic and transversely isotropic elastic layer bonded to a homogeneous, linearly elastic and transversely isotropic half-space is considered. Using the general solution to the governing integral equation of the axisymmetric contact problem for an isotropic elastic half-space, we derive exact equations for the contact force and the contact radius, which are then approximated under the assumption that the contact radius is sufficiently small compared to the thickness of the elastic layer. An asymptotic analysis of the resulting non-linear algebraic problem corresponding to the fourth-order asymptotic model is performed.
A special case of the indentation problem for a blunt punch of power-law profile is studied in detail. Approximate force-displacement relations are obtained in explicit form, which is most suited for development of indentation tests. 
\end{abstract}

\setcounter{equation}{0}

\section{Introduction}
\label{SecI}

In recent years, development of minimally invasive measurement techniques for determining the biomechanical properties of biological and artificial replacement tissues has become a problem of paramount importance in medicine. The mechanical properties of biomaterials can be evaluated by means of simple indentation techniques \cite{Fischer-Cripps2004}, which have been proved useful for both identification of mechanical properties of a soft biological tissue, like articular cartilage \cite{Hayes_et_al_1972}, and assessing its viability during arthroscopy \cite{Korhonen_et_al_2003}. Indentation tests are performed with various indenter geometries, and among the most frequently used types of indenters are axisymmetric (in particular, cylindrical, spherical, and conical) indenters. 

\begin{figure}[h!]
%\vskip-1.0cm    
    \centering
    \hbox{
    \includegraphics[scale=0.4]{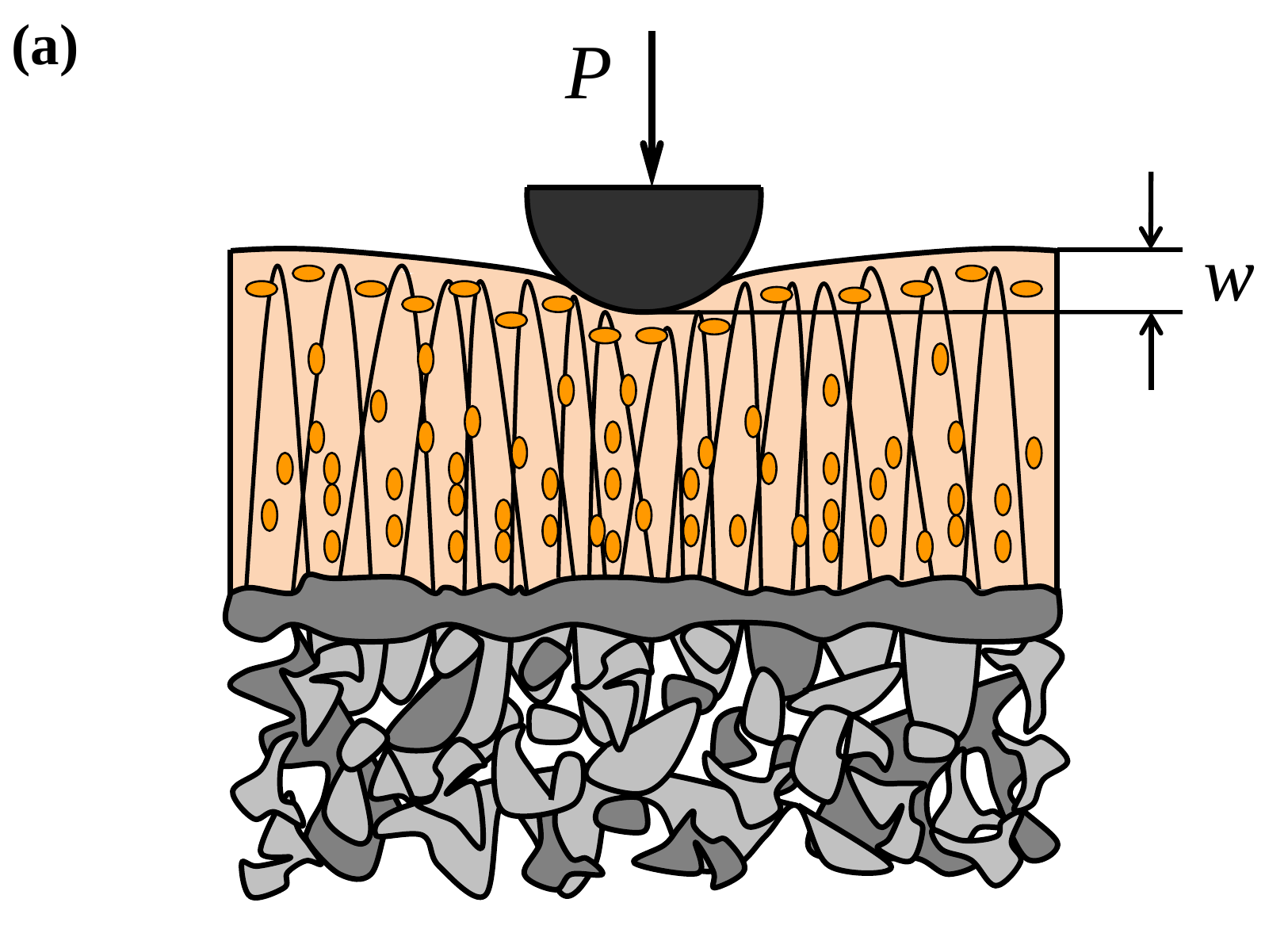}\hskip0.5cm
    \includegraphics[scale=0.4]{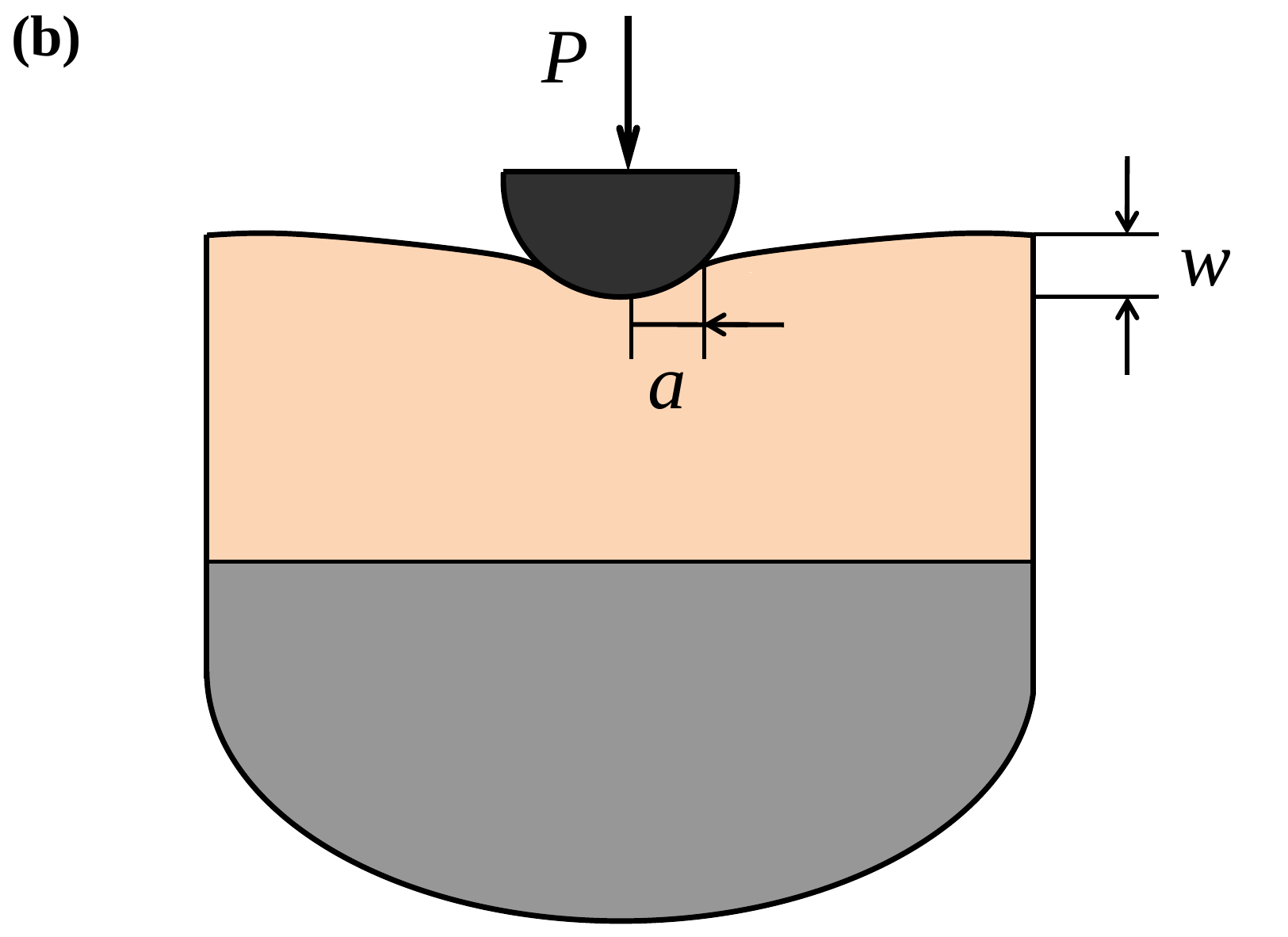}
    }
%\vskip-3.5cm    
    \caption{a) Schematics of indentation test for the articular cartilage layer attached to subchondral bone;
    b) Schematics of the indentation problem for an elastic layer bonded to an elastic half-space.}
%\vskip-1.0cm        
    \label{FigArbitray}
\end{figure}

Development of an indentation test for the articular cartilage layer requires solving the corresponding contact problem for an elastic layer (see Fig.~\ref{FigArbitray}a). The widely used mathematical model for indentation tests of articular cartilage developed by Hayes {\it et al.} \cite{Hayes_et_al_1972} is based on the analytical solution \cite{LebedevUfliand1958} of the frictionless contact problem for an elastic layer bonded to a rigid base. The Hayes {\it et al.} model allows to take into account the thickness effect, but, at the same time, it neglects the substrate effect. In order to take into account the influence of the substrate deformation, we consider the frictionless contact problem for an elastic layer bonded to an elastic half-space (see Fig.~\ref{FigArbitray}b). 

In the present paper, we consider the quasistatic unilateral contact problem for a rigid axisymmetric indenter pressed into a homogeneous, linearly elastic and transversely isotropic elastic layer bonded to a homogeneous, linearly elastic and transversely isotropic half-space. We focus on deriving asymptotic approximations for the force-displacement relation (that is the relation between the contact force $P$ and the indenter displacement $w$), which constitutes the basis for developing depth-sensing indentation tests. In the isotropic case, the substrate effect in the flat-ended and spherical indentation was studied recently in \cite{Argatov2010nano} based on the asymptotic solutions obtained in \cite{Vorovich_et_al_1974,Argatov2001finite}.
The main novelty of the present paper consists in considering the case of an arbitrary axisymmetric indenter. As a special case, we study the indentation problem for a blunt punch of power-law profile. 

The main difficulty for an analytical analysis of the axisymmetric unilateral contact problem is associated with the fact that the radius $a$ of the circular contact area $\omega$ is not known in advance and must be determined as a part of the solution. Apparently the first asymptotic analysis of the axisymmetric contact problem for an arbitrary axisymmetric indenter dates back to England, whose work \cite{England1962} was published in 1962. As a special case, England \cite{England1962} considered the indentation problem for a hemispherically ended punch and obtained asymptotic expansions in terms of powers of a small parameter $\mu=R/h$, where $R$ is the curvature radius of the hemispherical end, $h$ is the thickness of the elastic layer. 

Our asymptotic analysis is based on a small parameter $\varepsilon=a/h$ defined as the ratio of the contact radius to the layer thickness. It should be noted that in contrast to the method originally developed by England \cite{England1962}, we construct asymptotic expansions in terms of an unknown quantity in the same way as in the asymptotic method of large $\lambda$ developed by 
Vorovich {\it et al.} \cite{Vorovich_et_al_1974}. However, our approach differs from previous studies \cite{Vorovich_et_al_1974,AlexandrovPozharskii2001} in two ways. First, we do not construct an asymptotic expansion for the contact pressure, but for the aim of deriving the force-displacement relation, following \cite{Arg2004layer,ArgatovSabina2012}, we derive exact equations for the contact force and the contact radius and operate with integral characteristics of the contact pressure. Second, we perform asymptotic analysis of the resulting non-linear algebraic problem obtaining the approximate force-displacement relation in explicit form, which is most suited for development of indentation tests. 

\section{Indentation problem formulation}
\label{Sec1}

We consider the unilateral frictionless contact problem for an elastic layer bonded to an elastic half-space, which is made of a different material. Let the composite elastic medium occupying the half-space $x_3>0$ be indented by an arbitrary frictionless axisymmetric indenter. Assuming that the $x_3$-axis coincides with the axis of symmetry of the indenter, we will describe its surface by the equation 
\begin{equation}
x_3=\Phi(x_1,x_2),
\label{a1z(1.1)}
\end{equation}
where $\Phi(x_1,x_2)$ is the shape function whose dependence on the Cartesian coordinates $x_1$ and $x_2$ comes through a polar radius $r=\sqrt{x_1^2+x_2^2}$. Without loss of generality, we may assume that $\Phi(0,0)=0$.

Under the assumption that the indenter shape function $\Phi(x_1,x_2)$ is convex and the indenter is loaded with a normal force, $P$, which is directed along the positive $x_3$-axis, contact between the indenter and the surface of the elastic medium will be established over a circular region, $\omega$, say, of radius $a$. For a blunt indenter, the contact area $\omega$ is variable and depends on the contact force $P$.

Further, let $G_3(x_1,x_2,0)$ be the surface influence function for the elastic medium that gives the vertical surface displacement of the elastic layer under a unit point force applied at the origin of coordinates and directed along $x_3$-axis. Then, the unilateral contact problem under consideration, which involves an a priori unknown contact area $\omega$, reduces to that of solving the following integral equation:
\begin{equation}
\iint\limits_\omega G_3(x_1-y_1,x_2-y_2,0)p({\bf y})\,d{\bf y}=w-\Phi(x_1,x_2).
\label{a1z(1.2)}
\end{equation}
Here, $p(x_1,x_2)$ is the contact pressure, $w$ is the indenter displacement. 

It should be emphasized that Eq.~(\ref{a1z(1.2)}) implicitly contains an a priori unknown radius $a$ of the contact area $\omega$, which is to be determined as a part of solution from the following conditions:
\begin{equation}
p({\bf y})>0 \ \ {\rm for} \ \ 0\leq \vert{\bf y}\vert<a, \quad
p({\bf y})=0 \ \ {\rm for} \ \ \vert{\bf y}\vert=a.
\label{a1z(1.3)}
\end{equation}
Here, $\vert{\bf y}\vert=\sqrt{y_1^2+y_2^2}$ is a polar radius. 

By the equilibrium equation, the external contact force $P$ is related to the contact pressure density $p(x_1,x_2)$ as follows:
\begin{equation}
\iint\limits_\omega p({\bf y})\,d{\bf y}=P.
\label{a1z(1.4)}
\end{equation}

Due to the axisymmetry of the elastic medium, the surface influence function can be represented as 
\begin{equation}
G_3(y_1,y_2,0)=\frac{1}{2\pi\theta}\biggl(\frac{1}{\vert{\bf y}\vert}-
\frac{1}{h}\mathcal{F}\Bigl(\frac{\vert{\bf y}\vert}{h}\Bigr)\biggr),
\label{a1z(1.5)}
\end{equation}
where $h$ is the layer thickness, $\theta$ is an elastic constant related to the first layer, $\mathcal{F}(t)$ is a dimensionless function.

It can be shown \cite{Vorovich_et_al_1974} that the following integral representation holds:
\begin{equation}
\mathcal{F}(t)=\int\limits_0^\infty
\bigl[1-L(u)\bigr]J_0(u t)\,du.
\label{a1z(1.6)}
\end{equation}
Here, $J_0(x)$ is the Bessel function of the first kind, the kernel function $L(u)$ depends both on the material properties of the elastic layer and the elastic half-space.

According to \cite{Vorovich_et_al_1974}, there is a neighborhood of zero on which the function 
(\ref{a1z(1.6)}) is expanded into an absolutely convergent power series 
\begin{equation}
\mathcal{F}(t)=\sum_{m=0}^\infty a_m t^{2m}
\label{a1z(1.7)}
\end{equation}
with the coefficients 
\begin{equation}
a_m=\frac{(-1)^m}{[(2m)!!]^2}\int\limits_0^\infty
\bigl[1-L(u)\bigr]u^{2m}\,du.
\label{a1z(1.8)}
\end{equation}

Thus, the indentation problem is to find the force-displacement relationship, $P(w)$. In the present paper, we employ an asymptotic approach to derive approximate analytical formulas for $P(w)$ in the small-contact limit, that is under the assumption that the contact radius $a$ is small compared with the layer thickness $h$.

\section{Equation for the radius of the contact area}
\label{Sec2}

In the axisymmetric case, Eq.~(\ref{a1z(1.2)}) can be reduced to a simpler form (see, e.\,g.,
\cite{Vorovich_et_al_1974,AlexandrovPozharskii2001})
\begin{equation}
\frac{2}{\pi}\int_0^a 
{\bf K}\biggl(\frac{2\sqrt{\rho r}}{\rho+r}\biggr)\frac{p(\rho)\rho}{r+\rho}
\,d\rho=\theta\bigl(w-\Phi(r)\bigr)+\frac{1}{h}\int_0^a 
\mathcal{F}\Bigl(\frac{\rho}{h},\frac{r}{h}\Bigr)p(\rho)\rho\,d\rho,
\label{a1z(2.1)}
\end{equation}
where ${\bf K}(x)$ is the complete elliptic integral of the first kind, $\mathcal{F}(\sigma,\tau)$ is a dimensionless function given by
\begin{equation}
\mathcal{F}(\sigma,\tau)=\int_0^\infty 
\bigl[1-L(u)\bigr]J_0(\sigma u)J_0(\tau u)\,du.
\label{a1z(2.2)}
\end{equation}

In the neighborhood of the origin, the function $\mathcal{F}(\sigma,\tau)$ has the following absolutely convergent expansion \cite{AlexandrovPozharskii2001}:
\begin{equation}
\mathcal{F}(\sigma,\tau)=\sum_{i=0}^\infty\sum_{j=0}^\infty
b_{ij}\sigma^{2i}\tau^{2j},
\label{a1z(2.3)}
\end{equation}
$$
b_{ij}=\frac{(m!)^2}{(i!)^2(j!)^2}a_m \quad (m=i+j),
$$
where $a_m$ are given by (\ref{a1z(1.8)}).

Now, introducing an auxiliary function $v(r)$ by the formula 
\begin{equation}
\frac{1}{2\pi}v(r)=\theta \bigl(w-\Phi(r)\bigr)+
\frac{1}{h}\int\limits_0^a
\mathcal{F}\left(\frac{\rho}{h},\frac{r}{h}\right) p(\rho)\rho\,d\rho
\label{a1z(2.4)}
\end{equation}
and making use of the general solution of the axisymmetric contact problem for an elastic isotropic half-space \cite{Leonov1939,Schubert1942,Shtaerman1949}, we represent the solution to 
Eq.~(\ref{a1z(2.1)}) in the form
\begin{equation}
p(r)=\frac{\mathcal{V}(a)}{\pi\sqrt{a^2-r^2}}-\frac{1}{\pi}\int\limits_r^a
\frac{\mathcal{V}^{\prime}(s)\,ds}{\sqrt{s^2-r^2}},
\label{a1z(2.5)}
\end{equation}
\begin{equation}
\pi \mathcal{V}(r)=v(0)+r\int\limits_0^r
\frac{v^{\prime}(s)\,ds}{\sqrt{r^2-s^2}}. 
\label{a1z(2.6)}
\end{equation}

From the boundary condition (\ref{a1z(1.3)}) of vanishing the contact pressure at the edge contour of the contact area, it is necessarily follows that 
\begin{equation}
\mathcal{V}(a)=0.
\label{a1z(2.7)}
\end{equation}

Hence, in view of (\ref{a1z(2.7)}), formula (\ref{a1z(2.5)}) reduces to 
\begin{equation}
p(r)=-\frac{1}{\pi}\int\limits_r^a
\frac{\mathcal{V}^{\prime}(s)\,ds}{\sqrt{s^2-r^2}}.
\label{a1z(2.8)}
\end{equation}

Now, according to Eq.~(\ref{a1z(2.4)}) and the assumption $\Phi(0)=0$, we will have
\begin{equation}
\frac{1}{2\pi}v(0)=\theta w+ \frac{1}{h}\int\limits_0^a
\mathcal{F}\left(\frac{\rho}{h},0\right) p(\rho)\rho\,d\rho, 
\label{a1z(2.9)}
\end{equation}
\begin{equation}
\frac{1}{2\pi}v^{\prime}(r)=-\theta \Phi^{\prime}(r)+
\frac{1}{h^2}\int\limits_0^a \frac{\partial \mathcal{F}}{\partial\tau}
\left(\frac{\rho}{h},\frac{r}{h}\right)p(\rho)\rho\,d\rho.
\label{a1z(2.10)}
\end{equation}

The substitution of the expression (\ref{a1z(2.2)}) into Eqs.~(\ref{a1z(2.9)}), (\ref{a1z(2.10)}) and then into Eq.~(\ref{a1z(2.6)}) yields the following formulas \cite{Arg2004layer}:
\begin{equation}
\mathcal{V}(r)=2\theta w-2\theta r\int\limits_0^r \frac{\Phi^{\prime}(s)\,ds}{\sqrt{r^2-s^2}}
+\frac{2}{h}\int\limits_0^a S_0\left(\frac{\rho}{h},\frac{r}{h}\right)
p(\rho)\rho\,d\rho,
\label{a1z(2.11)}
\end{equation}
\begin{equation}
S_0(\sigma,\tau)=
\int\limits_0^{\infty}\bigl[1-L(u)\bigr]J_0(\sigma u)\cos\tau u\,du. 
\label{a1z(2.12)}
\end{equation}

Thus, Eq.~(\ref{a1z(2.7)}) for determining the contact radius $a$ takes the form 
\begin{equation}
w=a\int\limits_0^a \frac{\Phi^{\prime}(\rho)\,d\rho}{\sqrt{a^2-\rho^2}}-
\frac{1}{\theta h}\int\limits_0^a S_0\left(\frac{\rho}{h},\frac{a}{h}\right)
p(\rho)\rho\,d\rho.
\label{a1z(2.13)}
\end{equation}

It is clear that if an approximation for the contact pressure $p(r)$ is known, the second term of the right-hand side of Eq.~(\ref{a1z(2.13)}) will give a correction to the corresponding equation of the Galin\,---\,Sneddon theory of axisymmetric elastic contact \cite{Galin2008,Sneddon1965}.

\section{Equation for the contact force}
\label{Sec3}

Substituting the expression (\ref{a1z(2.8)}) for the contact pressure into the equilibrium equation (\ref{a1z(1.4)}), we will have
\begin{equation}
P=2\int\limits_0^a \mathcal{V}(s)\,ds.
\label{a1z(3.1)}
\end{equation}

Now, the substitution of (\ref{a1z(2.11)}) for $\mathcal{V}(s)$ into Eq.~(\ref{a1z(3.1)}) yields 
\begin{equation}
P=4\theta w a-4\theta\int\limits_0^a\Phi^{\prime}(\rho)\sqrt{a^2-\rho^2}d\rho
+4\int\limits_0^a
S_2\left(\frac{\rho}{h},\frac{a}{h}\right)p(\rho)\rho\,d\rho,
\label{a1z(3.2)}
\end{equation}
where
\begin{equation}
S_2(\sigma,\alpha)=\int\limits_0^{\infty}\bigl[1-L(u)\bigr]J_0(\sigma
u) \frac{\sin\alpha u}{u}\,du.
\label{a1z(3.3)}
\end{equation}

Finally, taking into account Eq.~(\ref{a1z(2.13)}), we eliminate the variable $w$ in Eq.~(\ref{a1z(3.2)}) as follows:
\begin{equation}
P=4\theta\int\limits_0^a\frac{\Phi^{\prime}(\rho)\rho^2 d\rho}{\sqrt{a^2-\rho^2}}
+4\int\limits_0^a
\tilde{S}_2\left(\frac{\rho}{h},\frac{a}{h}\right)p(\rho)\rho\,d\rho,
\label{a1z(3.4)}
\end{equation}
Here we introduced the notation 
\begin{equation}
\tilde{S}_2(\sigma,\alpha)=\int\limits_0^{\infty}\bigl[1-L(u)\bigr]
J_0(\sigma u) \biggl(\frac{\sin\alpha u}{u}-\alpha\cos \alpha u\biggr)du.
\label{a1z(3.5)}
\end{equation}

Again, if we know have an approximation for the contact pressure $p(r)$, then the second term of the right-hand side of Eq.~(\ref{a1z(3.4)}) will produce a correction to the corresponding equation of the Galin\,---\,Sneddon theory of axisymmetric elastic contact \cite{Galin2008,Sneddon1965}.

\section{The fourth-order asymptotic model}
\label{Sec4}

According to Eqs.~(\ref{a1z(2.12)}) and (\ref{a1z(3.5)}), the following series expansions hold:
\begin{equation}
S_0(\sigma,\alpha)=\sum_{i=0}^\infty\sum_{j=0}^\infty
b_{ij}^{(0)}\sigma^{2i}\alpha^{2j},
\label{a1z(4.1)}
\end{equation}
$$
b_{ij}^{(0)}=\frac{2^{2j}(m!)^2}{(i!)^2(2j)!}a_m \quad (m=i+j),
$$
\begin{equation}
\tilde{S}_2(\sigma,\alpha)=\alpha\sum_{i=0}^\infty\sum_{j=1}^\infty
\tilde{b}_{ij}^{(2)}\sigma^{2i}\alpha^{2j},
\label{a1z(4.2)}
\end{equation}
$$
\tilde{b}_{ij}^{(2)}=-\frac{2^{2j+1}j(m!)^2}{(i!)^2(2j+1)!}a_m \quad (m=i+j).
$$

Now, keeping terms of the infinite series (\ref{a1z(4.1)}) and (\ref{a1z(4.2)}) which contain only the coefficients $a_0$ and $a_1$, we get
\begin{equation}
S_0(\sigma,\alpha)\simeq a_0+a_1(\sigma^2+2\alpha^2),
\label{a1z(4.3)}
\end{equation}
\begin{equation}
\tilde{S}_2(\sigma,\alpha)\simeq -\frac{4}{3}a_1\alpha^3.
\label{a1z(4.4)}
\end{equation}

The substitution of the asymptotic expressions (\ref{a1z(4.3)}) and (\ref{a1z(4.4)}) into Eqs.~(\ref{a1z(2.13)}) and (\ref{a1z(3.4)}) results in the following approximate relations:
\begin{equation}
w\simeq a\int\limits_0^a \frac{\Phi^{\prime}(\rho)\,d\rho}{\sqrt{a^2-\rho^2}}-
\frac{1}{2\pi\theta h}\biggl\{P\Bigl(a_0+2a_1\frac{a^2}{h^2}\Bigr)
-P_2\frac{a_1}{h^2}\biggr\},
\label{a1z(4.5)}
\end{equation}
\begin{equation}
P\simeq 4\theta\int\limits_0^a\frac{\Phi^{\prime}(\rho)\rho^2}{\sqrt{a^2-\rho^2}}\,d\rho
-\frac{8a_1}{3\pi}\frac{a^3}{h^3}P.
\label{a1z(4.6)}
\end{equation}
Here, $P$ and $P_2$ are, respectively, the contact force and the polar moment of inertia of the contact pressure, i.\,e.,
$$
P=2\pi\int\limits_0^a p(\rho)\rho\,d\rho,
$$
\begin{equation}
P_2=2\pi\int\limits_0^a \rho^2 p(\rho)\rho\,d\rho.
\label{a1z(4.7)}
\end{equation}

Substituting the expression (\ref{a1z(2.8)}) for the contact pressure into Eq.~(\ref{a1z(4.7)}), we obtain 
\begin{equation}
P_2=4\int\limits_0^a \mathcal{V}(s)s^2\,ds.
\label{a1z(4.8)}
\end{equation}

Now, the substitution of (\ref{a1z(2.11)}) for $\mathcal{V}(s)$ into Eq.~(\ref{a1z(4.8)}) yields 
\begin{equation}
P_2=8\theta \Biggl(\frac{a^3}{3}w-\int\limits_0^a\Phi(\rho)
\frac{(2\rho^2-a^2)}{\sqrt{a^2-\rho^2}}\rho\,d\rho\Biggr)
+\frac{8a^3}{h}\int\limits_0^a
S_4\left(\frac{\rho}{h},\frac{a}{h}\right)p(\rho)\rho\,d\rho.
\label{a1z(4.8b)}
\end{equation}
Here we introduced the notation 
\begin{equation}
S_4(\sigma,\alpha)=\int\limits_0^{\infty}\bigl[1-L(u)\bigr]J_0(\sigma u) 
\biggl\{
\frac{\sin\alpha u}{\alpha u}-\frac{2}{(\alpha u)^2}
\Bigl(\frac{\sin\alpha u}{\alpha u}-\cos\alpha u\Bigr)\biggr\}du.
\label{a1z(4.9)}
\end{equation}

In view of the fact that the quantity $P_2$ enters Eq.~(\ref{a1z(4.5)}) with the coefficient $a_1$, Eq.~(\ref{a1z(4.8b)}) may be reduced to the following one being in the same range of accuracy with Eqs.~(\ref{a1z(4.5)}) and (\ref{a1z(4.6)}):
\begin{equation}
P_2\simeq 8\theta \Biggl(\frac{a^3}{3}w-\int\limits_0^a\Phi(\rho)
\frac{(2\rho^2-a^2)}{\sqrt{a^2-\rho^2}}\rho\,d\rho\Biggr)+\frac{4a^3 a_0}{3\pi h}P.
\label{a1z(4.10)}
\end{equation}
Note that in writing down Eq.~(\ref{a1z(4.10)}), it is assumed that 
\begin{equation}
S_4(\sigma,\alpha)\simeq \frac{1}{3}a_0.
\label{a1z(4.11)}
\end{equation}

Thus, the three relations (\ref{a1z(4.5)}), (\ref{a1z(4.6)}), and (\ref{a1z(4.10)}) form the resulting algebraic problem with respect to the four variables $w$, $P$, $P_2$, and $a$.

Eliminating the variable $P_2$ from Eq.~(\ref{a1z(4.5)}) by the substitution (\ref{a1z(4.8)}), we obtain 
\begin{eqnarray}
w\biggl(1+\frac{4a_1}{3\pi}\frac{a^3}{h^3}\biggr) & \simeq &
a\int\limits_0^a \frac{\Phi^{\prime}(\rho)\,d\rho}{\sqrt{a^2-\rho^2}}+
\frac{4a_1}{\pi h^3}\int\limits_0^a\Phi(\rho)
\frac{(2\rho^2-a^2)}{\sqrt{a^2-\rho^2}}\rho\,d\rho
\nonumber \\
{} & {} & {}-\frac{P}{2\pi\theta h}\biggl(a_0+2a_1\frac{a^2}{h^2}
+\frac{4a_0 a_1}{3\pi}\frac{a^3}{h^3}\biggr).
\label{a1z(4.12)}
\end{eqnarray}

On the other hand, Eq.~(\ref{a1z(4.6)}) yields 
\begin{equation}
P\biggl(1+\frac{8a_1}{3\pi}\frac{a^3}{h^3}\biggr)
\simeq 4\theta\int\limits_0^a \frac{\Phi^{\prime}(\rho)\rho^2 d\rho}{\sqrt{a^2-\rho^2}}.
\label{a1z(4.13)}
\end{equation}

With the same asymptotic accuracy as that of Eqs.~(\ref{a1z(4.5)}), (\ref{a1z(4.6)}), and (\ref{a1z(4.10)}), we obtain from Eq.~(\ref{a1z(4.13)}) the following approximate relation:
\begin{equation}
P\simeq 
\biggl(1-\frac{8a_1}{3\pi}\frac{a^3}{h^3}\biggr)
4\theta\int\limits_0^a \frac{\Phi^{\prime}(\rho)\rho^2 d\rho}{\sqrt{a^2-\rho^2}}.
\label{a1z(4.14)}
\end{equation}

Now, substituting the expression (\ref{a1z(4.14)}) into Eq.~(\ref{a1z(4.12)}), we arrive after some asymptotic simplifications at the following expansion:
\begin{eqnarray}
w & \simeq & \biggl(1-\frac{4a_1}{3\pi}\frac{a^3}{h^3}\biggr)
a\int\limits_0^a \frac{\Phi^{\prime}(\rho)\,d\rho}{\sqrt{a^2-\rho^2}}+
\frac{4a_1}{\pi h^3}\int\limits_0^a\Phi(\rho)
\frac{(2\rho^2-a^2)}{\sqrt{a^2-\rho^2}}\rho\,d\rho
\nonumber \\
{} & {} & {}-\frac{2}{\pi h}\biggl(a_0+2a_1\frac{a^2}{h^2}
-\frac{8a_0 a_1}{3\pi}\frac{a^3}{h^3}\biggr)
\int\limits_0^a \frac{\Phi^{\prime}(\rho)\rho^2 d\rho}{\sqrt{a^2-\rho^2}}.
\label{a1z(4.15)}
\end{eqnarray}

Thus, Eqs.~(\ref{a1z(4.14)}) and (\ref{a1z(4.15)}) define the force displacement relation in a parametric way. 

We note that the fourth-order asymptotics of the axisymmetric unilateral contact problem for an elastic layer was previously derived by Alexandrov and Pozharskii \cite{AlexandrovPozharskii2001} in a somewhat different way (see, \cite{AlexandrovPozharskii2001}, Ch.~1.3, Eqs.~(25) and (26)).
It can be shown that the two solutions are asymptotically equivalent up to terms of order 
$O(\varepsilon^5)$ as $\varepsilon\to 0$, where $\varepsilon=a/h$.

\section{The case of a blunt indenter}
\label{Sec5}

Let us consider indentation of the elastic medium by a blunt indenter with the shape function 
\begin{equation}
\Phi(r)=A r^\lambda,
\label{a1z(5.1)}
\end{equation}
where $1<\lambda$ is an arbitrary real number, $A$ is a constant having the dimension 
$[{\rm L}^{1-\lambda}]$ with ${\rm L}$ being the dimension of length. 

Taking into account the integral 
$$
\int\limits_0^1\frac{t^{\lambda+1}dt}{\sqrt{1-t^2}}=
\frac{2^{\lambda-2}\lambda}{(\lambda+1)}\Gamma\Bigl(\frac{\lambda}{2}\Bigr)^2
\Gamma(\lambda)^{-1},
$$
we reduce Eqs.~(\ref{a1z(4.14)}) and (\ref{a1z(4.15)}) to the following ones, respectively: 
\begin{equation}
P\simeq \theta A F_1(\lambda) a^{\lambda+1}
\Bigl(1-\varepsilon^3\frac{8a_1}{3\pi}\Bigr),
\label{a1z(5.2)}
\end{equation}
\begin{equation}
w\simeq A F_2(\lambda) a^\lambda
\biggl\{1-\varepsilon\frac{2\lambda a_0}{\pi(\lambda+1)}
-\varepsilon^3\frac{8\lambda(2\lambda+5)a_1}{3\pi(\lambda+3)(\lambda+1)}
+\varepsilon^4\frac{16\lambda a_0 a_1}{3\pi^2(\lambda+1)}\biggr\}.
\label{a1z(5.3)}
\end{equation}
Here we introduced the notation 
\begin{equation}
F_1(\lambda)=\frac{\lambda^2 2^\lambda}{(\lambda+1)}\Gamma\Bigl(\frac{\lambda}{2}\Bigr)^2
\Gamma(\lambda)^{-1},
\label{a1z(5.4)}
\end{equation}
\begin{equation}
F_2(\lambda)=\lambda 2^{\lambda-2}\Gamma\Bigl(\frac{\lambda}{2}\Bigr)^2
\Gamma(\lambda)^{-1},
\label{a1z(5.5)}
\end{equation}
\begin{equation}
\varepsilon=\frac{a}{h}.
\label{a1z(5.6)}
\end{equation}

Following \cite{Argatov2005} and using the asymptotic method based on Lagrange's formula for solving algebraic equations \cite{DeBruijn1958}, we invert Eq.~(\ref{a1z(5.3)}) to find the contact radius as a function of the indenter displacement in following form
\begin{equation}
\frac{a}{h}\simeq \varpi\bigl(1+B_1 \varpi+B_2 \varpi^2+B_3 \varpi^3+B_4 \varpi^4\bigr)
\label{a1z(5.7)}
\end{equation}
with the coefficients given by 
$$
B_1=\frac{2 a_0}{\pi(\lambda+1)},\quad
B_2=\frac{(\lambda+3)}{2(\lambda+1)^2}\Bigl(\frac{2 a_0}{\pi}\Bigr)^2,
$$
$$
B_3=\frac{(\lambda+4)(\lambda+2)}{3(\lambda+1)^3}\Bigl(\frac{2 a_0}{\pi}\Bigr)^3
+\frac{(2\lambda+5)}{(\lambda+1)(\lambda+3)}\frac{8 a_1}{3\pi},
$$
$$
B_4=\frac{(2\lambda+5)(\lambda+5)(3\lambda+5)(\lambda+3)}{24(\lambda+1)^4(\lambda+3)}
\Bigl(\frac{2 a_0}{\pi}\Bigr)^4
+\frac{(\lambda^2+11\lambda+22)}{(\lambda+1)^2(\lambda+3)}\frac{16 a_0 a_1}{3\pi^2}.
$$
Here we introduced the notation 
\begin{equation}
\varpi=\frac{1}{h}\Bigl(\frac{w}{A F_2(\lambda)}
\Bigr)^{\frac{1}{\lambda}}.
\label{a1z(5.8)}
\end{equation}

Now, substituting the expansion (\ref{a1z(5.7)}) into Eq.~(\ref{a1z(5.2)}), we obtain 
\begin{equation}
P\simeq \theta A^{-\frac{1}{\lambda}}F_3(\lambda) w^{\frac{\lambda+1}{\lambda}}
\bigl(1+C_1 \varpi+C_2 \varpi^2+C_3 \varpi^3+C_4 \varpi^4\bigr),
\label{a1z(5.9)}
\end{equation}
where 
\begin{equation}
F_3(\lambda)=\lambda^{\frac{\lambda-1}{\lambda}}
2^{\frac{\lambda+2}{\lambda}}(\lambda+1)^{-1}
\Gamma\Bigl(\frac{\lambda}{2}\Bigr)^{-\frac{2}{\lambda}}
\Gamma(\lambda)^{\frac{1}{\lambda}},
\label{a1z(5.10)}
\end{equation}
$$
C_1=\frac{2 a_0}{\pi},\quad
C_2=\frac{(2\lambda+3)}{2(\lambda+1)}\Bigl(\frac{2 a_0}{\pi}\Bigr)^2,
$$
$$
C_3=\frac{(\lambda+2)(3\lambda+4)}{3(\lambda+1)^2}\Bigl(\frac{2 a_0}{\pi}\Bigr)^3
+\frac{(\lambda+2)}{(\lambda+3)}\frac{8 a_1}{3\pi},
$$
\begin{equation}
C_4=\frac{(2\lambda+5)(3\lambda+5)(4\lambda+5)}{24(\lambda+1)^3}
\Bigl(\frac{2 a_0}{\pi}\Bigr)^4
+\frac{(2\lambda+5)(\lambda+2)}{(\lambda+1)(\lambda+3)}\frac{16 a_0 a_1}{3\pi^2}.
\label{a1z(5.9a)}
\end{equation}

In view of (\ref{a1z(5.8)}), Eq.~(\ref{a1z(5.9)}) represents the load-displacement relationship for a blunt indenter with the shape function (\ref{a1z(5.1)}). Note that Eqs.~(\ref{a1z(5.7)}) and (\ref{a1z(5.9)}) generalize the second-order asymptotic model developed in \cite{Argatov2011}.

\subsection{Example. Paraboloid indenter}
\label{Sec5.1}

In the case $\lambda=2$ and $A=1/(2R)$, Eq.~(\ref{a1z(5.1)}) takes the form 
\begin{equation}
\Phi(r)=\frac{r^2}{2R},
\label{a1z(5.11)}
\end{equation}
where $R$ is the radius of curvature of the indenter�s surface at its vertex. 

According to Eqs.~(\ref{a1z(5.4)}), (\ref{a1z(5.5)}), (\ref{a1z(5.8)}), and (\ref{a1z(5.10)}), we will have
$$
F_1(2)=\frac{16}{3}, \quad
F_2(2)=2, \quad
F_3(2)=\frac{4\sqrt{2}}{3},
$$
$$
B_1=\frac{2a_0}{3\pi}, \quad
B_2=\frac{10a_0^2}{9\pi^2}, \quad
B_3=\frac{64a_0^3}{27\pi^3}+\frac{8a_1}{5\pi}, \quad
B_4=\frac{154a_0^4}{27\pi^4}+\frac{256 a_0 a_1}{45\pi^2},
$$
$$
\varpi=\frac{\sqrt{Rw}}{h},
$$
$$
C_1=\frac{2a_0}{\pi}, \quad
C_2=\frac{14a_0^2}{3\pi^2}, \quad
C_3=\frac{320 a_0^3}{27\pi^3}+\frac{32 a_1}{15\pi}, \quad
C_4=\frac{286 a_0^4}{9\pi^4}+\frac{64 a_0 a_1}{5\pi^2}.
$$
The obtained results are in agreement with those in \cite{Argatov2011}.

\subsection{Example. Cone indenter}
\label{Sec5.2}

For a cone indenter, we have $\lambda=1$ and $A=\tan\gamma$, whereas Eq.~(\ref{a1z(5.1)}) takes the form 
$$
\Phi(r)=r\tan\gamma,
$$
where $\gamma$ is the angle between the contact surface and the side surface of the cone. 

According to Eqs.~(\ref{a1z(5.4)}), (\ref{a1z(5.5)}), (\ref{a1z(5.8)}), and (\ref{a1z(5.10)}), now we have
$$
F_1(1)=\pi, \quad
F_2(1)=\frac{\pi}{2}, \quad
F_3(1)=\frac{4}{\pi},
$$
$$
B_1=\frac{a_0}{\pi}, \quad
B_2=\frac{2 a_0^2}{\pi^2}, \quad
B_3=\frac{5 a_0^3}{\pi^3}+\frac{7 a_1}{3\pi}, \quad
B_4=\frac{14 a_0^4}{\pi^4}+\frac{34 a_0 a_1}{3\pi^2},
$$
$$
\varpi=\frac{2w \cot\gamma}{\pi h},
$$
$$
C_1=\frac{2a_0}{\pi}, \quad
C_2=\frac{5 a_0^2}{\pi^2}, \quad
C_3=\frac{14 a_0^3}{\pi^3}+\frac{2 a_1}{\pi}, \quad
C_4=\frac{42 a_0^4}{\pi^4}+\frac{14 a_0 a_1}{\pi^2}.
$$

Substituting the above coefficients into Eqs.~(\ref{a1z(5.7)}) and (\ref{a1z(5.9)}), we arrive at the following approximate relations:
\begin{eqnarray}
a & \simeq & \frac{2\cot\gamma}{\pi}w \biggl\{
1+\frac{2a_0}{\pi^2}\cot\gamma\frac{w}{h}+
\frac{8a_0^2}{\pi^4}\cot^2\gamma\frac{w^2}{h^2}
\nonumber \\
{} & {} & {}+\biggl(\frac{40 a_0^3}{\pi^6}+\frac{56 a_1}{3\pi^4}\biggr)\cot^3\gamma\frac{w^3}{h^3}
+\biggl(\frac{224 a_0^4}{\pi^8}+\frac{544 a_0 a_1}{3\pi^6}\biggr)\cot^4\gamma\frac{w^4}{h^4}
\biggr\},
\label{a1z(5.12)}
\end{eqnarray}
\begin{eqnarray}
P & \simeq & \frac{4\theta\cot\gamma}{\pi} w^2 \biggl\{
1+\frac{4a_0}{\pi^2}\cot\gamma\frac{w}{h}+
\frac{20 a_0^2}{\pi^4}\cot^2\gamma\frac{w^2}{h^2}
\nonumber \\
{} & {} & {}+\biggl(\frac{112 a_0^3}{\pi^6}+\frac{16 a_1}{\pi^4}\biggr)\cot^3\gamma\frac{w^3}{h^3}
+\biggl(\frac{672 a_0^4}{\pi^8}+\frac{224 a_0 a_1}{\pi^6}\biggr)\cot^4\gamma\frac{w^4}{h^4}
\biggr\}.
\label{a1z(5.13)}
\end{eqnarray}

Equations (\ref{a1z(5.12)}) and (\ref{a1z(5.13)}) generalize the second-order asymptotic model developed in \cite{Argatov2011}. We note the following misprints in \cite{Argatov2011}: missing factors $\cot\gamma$ and $\cot^2\gamma$ in Eq.~(31) and an extra factor $\vartheta$ in Eq.~(32).

\section{The case of a hemispherically-ended indenter}
\label{Sec6}

For a hemispherical indenter, we will have
\begin{equation}
\Phi(r)=R-\sqrt{R^2-r^2},
\label{a1z(6.1)}
\end{equation}
where $R$ is the curvature radius of the hemispherical end.

In what follows, we make use of the identities 
$$
\int\limits_0^a\frac{\Phi^\prime(\rho)\,d\rho}{\sqrt{a^2-\rho^2}}=
\frac{1}{2}\ln\biggl(\frac{1+\alpha}{1-\alpha}\biggr),
$$
$$
\int\limits_0^a\frac{\Phi^\prime(\rho)\rho^2 d\rho}{\sqrt{a^2-\rho^2}}=R^2
\biggl\{-\frac{\alpha}{2}+\frac{1}{4}(1+\alpha^2)
\ln\biggl(\frac{1+\alpha}{1-\alpha}\biggr)\biggr\},
$$
$$
\int\limits_0^a \Phi(\rho)\frac{(2\rho^2-a^2)}{\sqrt{a^2-\rho^2}}
\rho\,d\rho=R^4
\biggl\{\frac{\alpha}{4}+\frac{\alpha^3}{12}-\frac{1}{8}(1-\alpha^4)
\ln\biggl(\frac{1+\alpha}{1-\alpha}\biggr)\biggr\}.
$$
Here we introduced the notation 
\begin{equation}
\alpha=\frac{a}{R}.
\label{a1z(6.2)}
\end{equation}

Thus, Eqs.~(\ref{a1z(4.14)}) and (\ref{a1z(4.15)}) can be represented as follows:
\begin{equation}
P\simeq \theta R^2\biggl((1+\alpha^2)\ln\biggl(\frac{1+\alpha}{1-\alpha}\biggr)-2\alpha\biggr)
\biggl(1-\mu^3\alpha^3\frac{8a_1}{3\pi}\biggr),
\label{a1z(6.3)}
\end{equation}
\begin{eqnarray}
\frac{w}{R} & \simeq & \frac{\alpha}{2}\ln\biggl(\frac{1+\alpha}{1-\alpha}\biggr)
-\mu\frac{a_0}{2\pi}\biggl((1+\alpha^2)\ln\biggl(\frac{1+\alpha}{1-\alpha}\biggr)-2\alpha\biggr)
\nonumber \\
{} & {} & {}+\mu^3\frac{a_1}{6\pi}\biggl(2\alpha(7\alpha^2+3)-(7\alpha^4+6\alpha^2+3)
\ln\biggl(\frac{1+\alpha}{1-\alpha}\biggr)\biggr)
\nonumber \\
{} & {} & {}+\mu^4\frac{4a_0 a_1}{3\pi^2}\alpha^3
\biggl((1+\alpha^2)\ln\biggl(\frac{1+\alpha}{1-\alpha}\biggr)-2\alpha\biggr).
\label{a1z(6.4)}
\end{eqnarray}
Here we introduced the notation 
\begin{equation}
\mu=\frac{R}{h}.
\label{a1z(6.5)}
\end{equation}

Now, following England \cite{England1962}, we construct the approximate solution for the contact radius in the form of a power expansion with respect to the small parameter $\mu$ as
\begin{equation}
\frac{a}{R}\simeq \alpha_0+\mu\alpha_1+\mu^2\alpha_2+\mu^3\alpha_3+\mu^4\alpha_4 .
\label{a1z(6.6)}
\end{equation}
The coefficients of the expansion (\ref{a1z(6.6)}) are determined according to Eq.~(\ref{a1z(6.3)}). By a perturbation method, we get
\begin{equation}
\alpha_1=0,\quad \alpha_2=0,\quad \alpha_4=0,
\label{a1z(6.7)}
\end{equation}
\begin{equation}
\alpha_3=\frac{4a_1}{3\pi}\frac{\alpha_0^2(1-\alpha_0^2)\bigl[
(1+\alpha_0^2)\ln\bigl(\frac{1+\alpha_0}{1-\alpha_0}\bigr)-2\alpha_0\bigr]
}{
(1-\alpha_0^2)\ln\bigl(\frac{1+\alpha_0}{1-\alpha_0}\bigr)+2\alpha_0},
\label{a1z(6.8)}
\end{equation}
while $\alpha_0$ is a unique root of the equation 
\begin{equation}
\frac{P}{\theta R^2}=(1+\alpha_0^2)\ln\biggl(\frac{1+\alpha_0}{1-\alpha_0}\biggr)-2\alpha_0.
\label{a1z(6.9)}
\end{equation}

In view of (\ref{a1z(6.7)}), it is readily seen that the expansion (\ref{a1z(6.6)}) simplifies to
\begin{equation}
\frac{a}{R}= \alpha_0+\mu^3\alpha_3+O(\mu^5).
\label{a1z(6.10)}
\end{equation}

The substitution of (\ref{a1z(6.10)}) into Eq.~(\ref{a1z(6.4)}) yields the following result after expanding up to the fourth order in $\mu$:
\begin{eqnarray}
\frac{w}{R} & = & \frac{\alpha_0}{2}\ln\biggl(\frac{1+\alpha_0}{1-\alpha_0}\biggr)
-\mu\frac{a_0}{2\pi}\biggl((1+\alpha_0^2)\ln\biggl(\frac{1+\alpha_0}{1-\alpha_0}\biggr)-2\alpha_0\biggr)
\nonumber \\
{} & {} & {}-\mu^3\frac{a_1}{6\pi}\biggl((3\alpha_0^4+2\alpha_0^2+3)\ln\biggl(\frac{1+\alpha_0}{1-\alpha_0}\biggr)-6\alpha_0(1+\alpha_0^2)\biggr)+O(\mu^5).
\label{a1z(6.11)}
\end{eqnarray}

Formulas (\ref{a1z(6.6)})\,--\,(\ref{a1z(6.11)}) are in complete agreement with the corresponding formulas from \cite{England1962}, provided the following relations hold:
$$
\biggl(\frac{k_1}{1+k_1}-\frac{k_2}{1+k_2}\biggr)\frac{1}{\gamma_1-\gamma_2}=\frac{1}{\theta},
\quad K_0=-\frac{2a_0}{\pi}, \quad K_1=-\frac{4a_1}{\pi}.
$$
Here, $k_1$, $k_2$, $\gamma_1$, $\gamma_2$ are elastic constants used in \cite{England1962},  $K_0$ and $K_1$ are the corresponding asymptotic constants.

\section{Indentation scaling factor for a blunt indenter}
\label{Sec7}

Let us rewrite Eq.~(\ref{a1z(5.9)}) in the form
\begin{equation}
P=\theta A^{-\frac{1}{\lambda}}F_3(\lambda) w^{\frac{\lambda+1}{\lambda}}
f_\lambda(\varpi),
\label{a1z(7.1)}
\end{equation}
introducing the indentation scaling factor $f_\lambda(\varpi)$. According to (\ref{a1z(5.9)}), we have
\begin{equation}
f_\lambda(\varpi)\simeq 
1+C_1 \varpi+C_2 \varpi^2+C_3 \varpi^3+C_4 \varpi^4,
\label{a1z(7.2)}
\end{equation}
where the coefficients $C_n$ are given by (\ref{a1z(5.9a)}).

Observe that the factor $f_\lambda(\varpi)$ depends on the dimensionless indentation parameter $\varpi$ defined by Eq.~(\ref{a1z(5.8)}). In view of Eqs.~(\ref{a1z(5.3)}) and (\ref{a1z(5.8)}), we have
\begin{eqnarray}
\varpi & \simeq & \varepsilon\biggl\{
1-\varepsilon\frac{2a_0}{\pi(\lambda+1)}
-\varepsilon^2\frac{(\lambda-1)}{2(\lambda+1)^2}\Bigl(\frac{2a_0}{\pi}\Bigr)^2
\nonumber \\
{} & {} & {}-\varepsilon^3\biggl(
\frac{(\lambda-1)(2\lambda-1)}{6(\lambda+1)^3}\Bigl(\frac{2a_0}{\pi}\Bigr)^3
+\frac{8 a_1(2\lambda+5)}{3\pi(\lambda+3)(\lambda+1)}\biggr)
\nonumber \\
{} & {} & {}-\varepsilon^4\biggl(
\frac{(\lambda-1)(2\lambda-1)(3\lambda-1)}{24(\lambda+1)^4}\Bigl(\frac{2a_0}{\pi}\Bigr)^4
+\frac{16 a_0 a_1}{3\pi^2}\frac{(\lambda^2-\lambda-8)}{(\lambda+1)^2(\lambda+3)}
\biggr)\biggr\}.
\label{a1z(7.3)}
\end{eqnarray}

Now, substituting the expansion (\ref{a1z(7.3)}) into (\ref{a1z(7.2)}), we express the indentation factor in terms of the relative contact radius by introducing a new notation 
\begin{equation}
\kappa_\lambda(\varepsilon)=f_\lambda(\varpi(\varepsilon)).
\label{a1z(7.4)}
\end{equation}

Thus, from (\ref{a1z(7.2)}) and (\ref{a1z(7.3)}), it follows that 
\begin{eqnarray}
\kappa_\lambda(\varepsilon) & \simeq & 
1+\varepsilon\frac{2a_0}{\pi}
+\varepsilon^2\frac{(2\lambda+1)}{2(\lambda+1)}\Bigl(\frac{2a_0}{\pi}\Bigr)^2
\nonumber \\
{} & {} & {}+\varepsilon^3\biggl(
\frac{(2\lambda+1)(3\lambda+1)}{6(\lambda+1)^2}\Bigl(\frac{2a_0}{\pi}\Bigr)^3
+\frac{8 a_1(\lambda-2)}{3\pi(\lambda+3)}\biggr)
\nonumber \\
{} & {} & {}+\varepsilon^4\biggl(
\frac{(2\lambda+1)(3\lambda+1)(4\lambda+1)}{24(\lambda+1)^3}\Bigl(\frac{2a_0}{\pi}\Bigr)^4
+\frac{16 a_0 a_1}{3\pi^2}\frac{(2\lambda^2+4\lambda-1)}{(\lambda+1)(\lambda+3)}
\biggr).
\label{a1z(7.5)}
\end{eqnarray}

Substituting the value $\lambda=2$ into Eq.~(\ref{a1z(7.5)}), we readily get
\begin{eqnarray}
\kappa_2(\varepsilon) & \simeq & 
1+\varepsilon\frac{2a_0}{\pi}
+\varepsilon^2\frac{10 a_0^2}{3\pi^2}
+\varepsilon^3\biggl(
\frac{140 a_0^3}{27\pi^3}+\frac{32 a_1}{15\pi}\biggr)
\nonumber \\
{} & {} & {}+\varepsilon^4\biggl(
\frac{70 a_0^4}{9\pi^4}
+\frac{16 a_0 a_1}{3\pi^2}\biggr).
\label{a1z(7.6)}
\end{eqnarray}

In order to check the asymptotic formula (\ref{a1z(7.5)}), we pass to the limit as $\lambda\to+\infty$, obtaining the following result:
\begin{eqnarray}
\kappa_\infty(\varepsilon) & \simeq & 
1+\varepsilon\frac{2a_0}{\pi}
+\varepsilon^2\Bigl(\frac{2a_0}{\pi}\Bigr)^2
+\varepsilon^3\biggl(
\Bigl(\frac{2a_0}{\pi}\Bigr)^3+\frac{8 a_1}{3\pi}\biggr)
\nonumber \\
{} & {} & {}+\varepsilon^4\biggl(
\Bigl(\frac{2a_0}{\pi}\Bigr)^4
+\frac{32 a_0 a_1}{3\pi^2}\biggr).
\label{a1z(7.7)}
\end{eqnarray}

The case $\lambda=+\infty$ corresponds to a cylindrical indenter. The asymptotic expansion (\ref{a1z(7.7)}) is in complete agreement with the fourth-order asymptotics derived previously in 
\cite{Vorovich_et_al_1974,Argatov2002CLC}.

\section{Displacement-force relationship for a blunt indenter}
\label{Sec8}

According to the notation (\ref{a1z(5.8)}), the force-displacement relation (\ref{a1z(5.9)}) can be rewritten as
\begin{equation}
P\simeq \theta A F_1(\lambda) h^{\lambda+1} \varpi^{\lambda+1}
\bigl(1+C_1 \varpi+C_2 \varpi^2+C_3 \varpi^3+C_4 \varpi^4\bigr),
\label{a1z(8.1)}
\end{equation}
where the factor $F_1(\lambda)$ is given by (\ref{a1z(5.4)}).

Inverting the relation (\ref{a1z(8.1)}), we obtain 
\begin{equation}
\varpi\simeq \tilde{P}
\bigl(1+D_1 \tilde{P}+D_2 \tilde{P}^2+D_3 \tilde{P}^3+D_4 \tilde{P}^4\bigr),
\label{a1z(8.2)}
\end{equation}
where we introduced the notation 
\begin{equation}
\tilde{P}=\bigl(\theta A F_1(\lambda)\bigr)^{-\frac{1}{\lambda+1}}
h^{-1} P^{\frac{1}{\lambda+1}}.
\label{a1z(8.3)}
\end{equation}

Now, in view of (\ref{a1z(5.8)}), from Eq.~(\ref{a1z(8.2)}) it follows that 
\begin{equation}
w\simeq A^{\frac{1}{\lambda+1}}\Bigl(
\frac{P}{\theta F_3(\lambda)}\Bigr)^{\frac{\lambda}{\lambda+1}}
\biggl(1-\frac{2\lambda a_0}{\pi(\lambda+1)}\tilde{P}
-\frac{8\lambda(\lambda+2) a_1}{3\pi(\lambda+1)(\lambda+3)}\tilde{P}^3\biggr).
\label{a1z(8.4)}
\end{equation}

Formula (\ref{a1z(8.4)}) represents the sought-for relation between the indenter displacement $w$ and the contact force $P$. 

Substituting now the expansion (\ref{a1z(8.2)}) into Eq.~(\ref{a1z(5.7)}), we derive the corresponding relation between the contact radius $a$ and the contact force $P$ in the simple form
\begin{equation}
\frac{a}{h}\simeq \tilde{P}\biggl(1+\frac{8 a_1}{3\pi(\lambda+1)}\tilde{P}^3\biggr).
\label{a1z(8.5)}
\end{equation}

Observe that Eq.~(\ref{a1z(8.5)}) can be directly obtained from Eq.~(\ref{a1z(5.2)}).

\section{Incremental indentation stiffness for a blunt indenter}
\label{Sec9}

During the depth-sensing indentation, the incremental indentation stiffness can be evaluated according to the relation 
\begin{equation}
\frac{dP}{dw}=\frac{\frac{dP}{da}}{\frac{dw}{da}}.
\label{a1z(9.1)}
\end{equation}

Substituting the expressions (\ref{a1z(5.2)}) and (\ref{a1z(5.3)}) into Eq.~(\ref{a1z(9.1)}), we arrive at the approximate relation
\begin{equation}
\frac{dP}{dw}\simeq 4\theta a\frac{1-\frac{8 a_1}{3\pi}\frac{(\lambda+4)}{(\lambda+1)}\varepsilon^3}{
1-\frac{2a_0}{\pi}\varepsilon
-\frac{8 a_1}{3\pi}\frac{(2\lambda+5)}{(\lambda+1)}\varepsilon^3
+\frac{16 a_0 a_1}{3\pi^2}\frac{(\lambda+4)}{(\lambda+1)}\varepsilon^4}.
\label{a1z(9.2)}
\end{equation}
Note that in deriving formula (\ref{a1z(9.2)}), the following relation was taken into account (see Eqs.~(\ref{a1z(5.4)}) and (\ref{a1z(5.5)})):
$$
\frac{(\lambda+1)}{\lambda}\frac{F_1(\lambda)}{F_2(\lambda)}=4.
$$

Now, expanding the right-hand side of Eq.~(\ref{a1z(9.2)}) into a power series in $\varepsilon$, we obtain
\begin{eqnarray}
\frac{dP}{dw} & \simeq & 4\theta a\biggl\{
1+\varepsilon\frac{2a_0}{\pi}
+\varepsilon^2\Bigl(\frac{2a_0}{\pi}\Bigr)^2
+\varepsilon^3\biggl(
\Bigl(\frac{2a_0}{\pi}\Bigr)^3+\frac{8 a_1}{3\pi}\biggr)
\nonumber \\
{} & {} & {}+\varepsilon^4\biggl(
\Bigl(\frac{2a_0}{\pi}\Bigr)^4
+\frac{32 a_0 a_1}{3\pi^2}\biggr)\biggr\}.
\label{a1z(9.3)}
\end{eqnarray}

It is interesting to observe that the coefficients of the asymptotic expansion (\ref{a1z(9.3)}) do not depend on the parameter $\lambda$, which describes the indenter shape. 

On the other hand, the incremental indentation stiffness for a blunt indenter can be evaluated by differentiating the force-displacement relationship (\ref{a1z(8.1)}) with respect to the indenter displacement. Thus, by taking into account the formula 
$$
\frac{d\varpi}{dw}=\frac{\varpi^{1-\lambda}}{\lambda F_2(\lambda)Ah},
$$
we get the relation
\begin{eqnarray}
\frac{dP}{dw} & \simeq & \theta \frac{(\lambda+1)}{\lambda}F_3(\lambda)
F_2(\lambda)^{\frac{1}{\lambda}}h\varpi\biggl\{
1+\varpi C_1\frac{(\lambda+2)}{\lambda+1}
+\varpi^2 C_2\frac{(\lambda+3)}{\lambda+1}
\nonumber \\
{} & {} & {}+\varpi^3 C_3\frac{(\lambda+4)}{\lambda+1}
+\varpi^4 C_4\frac{(\lambda+5)}{\lambda+1}\biggr\},
\label{a1z(9.4)}
\end{eqnarray}
where the coefficients $C_n$ are given by (\ref{a1z(5.9a)}).

Now, substituting the expression (\ref{a1z(7.3)}) for $\varpi$ into Eq.~(\ref{a1z(9.4)}) and expanding the result in terms of powers of the parameter $\varepsilon$, we again arrive at Eq.~(\ref{a1z(9.3)}) by noting that 
$$
\frac{(\lambda+1)}{\lambda}F_3(\lambda)F_2(\lambda)^{\frac{1}{\lambda}}=4.
$$

Finally, by comparing Eqs.~(\ref{a1z(9.3)}) and (\ref{a1z(7.7)}), it is readily seen that Eq.~(\ref{a1z(9.3)}) can be rewritten in the following form:
\begin{equation}
\frac{dP}{dw}=4\theta a\kappa_\infty(\varepsilon).
\label{a1z(9.5)}
\end{equation}
Here, $\kappa_\infty(\varepsilon)$ is the indentation scaling factor for a cylindrical indenter.

Observe that Eq.~(\ref{a1z(9.5)}) can be regarded as a generalization of Barber's theorem \cite{Barber2003} for the incremental indentation stiffness established in the case of an isotropic elastic half-space.

Note also that in context of indentation testing, Eq.~(\ref{a1z(9.5)}) should be rewritten as
\begin{equation}
\frac{dP}{dw}=4\theta \sqrt{\frac{A}{\pi}}\kappa_\infty(\varepsilon),
\label{a1z(9.6)}
\end{equation}
where $A$ is the current area of contact. 

In contrast to the Bulychev\,---\,Alekhin\,---\,Shorshorov (BASh) relation established in \cite{Bulychev_et_al_1976,Pharr_et_al_1992,BorodichKeer2004a} in the case of an isotropic elastic half-space, Eq.~(\ref{a1z(9.6)}) contains the extra factor $\kappa_\infty(\varepsilon)$, which accounts for both the thickness and substrate effects. Note that in the case of isotropy, we have
$\theta=E^{(1)}/\bigl[2\bigl(1-\bigl(\nu^{(1)}\bigr)^2\bigr)\bigr]$, where $E^{(1)}$ and $\nu^{(1)}$ are Young's modulus and Poisson's ratio of the layer material.

\section{Conclusion}
\label{SecC}

An asymptotic analysis of the substrate effect has been performed for deriving approximate solutions of the quasistatic indentation problem for an arbitrary frictionless indenter. Explicit asymptotic formulas are derived in the framework of the fourth-order asymptotic model, which has been previously shown to be effective for describing the thickness effect for a spherical indenter on an elastic layer. The results obtained can be applied for development of indentation tests.

\bigskip

{\bf Acknowledgements}. 
The support of Conacyt project number 129658, 
Instituto de Ciencia de Materiales, Estudios de Posgrado en Matem\'aticas, and Fenomec, UNAM is gratefully acknowledged. 
The research has been carried out during the Marie Curie Fellowship
of Dr.~I.\,I.~Argatov at Aberystwyth University supported by the European Union Seventh
Framework Programme under contract number PIIF-GA-2009-253055.
The authors thank Ms.~Ana P\'erez Arteaga and Mr.~Ramiro Ch\'avez Tovar for computational support.

\section*{Appendix: Kernel function for a transversely isotropic layer bonded to a transversely isotropic half-space}
\label{SecA}

The constitutive relationship for a transversely isotropic material referred to the Cartesian coordinates 
$(x_1,x_2,x_3)$ with the $Ox_1 x_2$ plane coinciding with the plane of elastic symmetry can be written in the matrix form as follows \cite{Elliott1948}:
$$
\left(
\begin{array}{c}
\sigma_{11} \\
\sigma_{22} \\
\sigma_{33} \\
\sigma_{13} \\
\sigma_{23} \\
\sigma_{12}
\end{array}
\right)=
\left[
\begin{array}{cccccc}
A_{11} & A_{12} & A_{13} & 0 & 0 & 0 \\
A_{12} & A_{11} & A_{13} & 0 & 0 & 0 \\
A_{13} & A_{13} & A_{33} & 0 & 0 & 0 \\
0 & 0 & 0 & 2A_{44} & 0 & 0 \\
0 & 0 & 0 & 0 & 2A_{44} & 0 \\
0 & 0 & 0 & 0 & 0 & A_{11}-A_{12}
\end{array}
\right]
\left(
\begin{array}{c}
\varepsilon_{11} \\
\varepsilon_{22} \\
\varepsilon_{33} \\
\varepsilon_{13} \\
\varepsilon_{23} \\
\varepsilon_{12}
\end{array}
\right).
$$
For a transversely isotropic material, only five independent elastic constants are needed to describe its deformational behavior. The elastic moduli $A_{11}$, $A_{12}$, $A_{13}$, $A_{33}$, and $A_{44}$ can be expressed in terms of the engineering elastic constants as follows \cite{LiaoWang1998}:
$$
A_{11}=\frac{E\Bigl(1-\frac{E}{E^\prime}\,\nu^{\prime 2}\Bigr)}{
(1+\nu)\Bigl(1-\nu-\frac{2E}{E^\prime}\,\nu^{\prime 2}\Bigr)}, \quad
A_{12}=\frac{E\Bigl(\nu+\frac{E}{E^\prime}\,\nu^{\prime 2}\Bigr)}{
(1+\nu)\Bigl(1-\nu-\frac{2E}{E^\prime}\,\nu^{\prime 2}\Bigr)},
$$
$$
A_{13}=\frac{E\nu^\prime}{1-\nu-\frac{2E}{E^\prime}\,\nu^{\prime 2}},\quad
A_{33}=\frac{E^\prime(1-\nu)}{
1-\nu-\frac{2E}{E^\prime}\,\nu^{\prime 2}}, \quad
A_{44}=G^\prime.
$$
Here, $E$ and $E^\prime$ are Young's moduli in the plane of transverse isotropy and in the direction normal to it, respectively, $\nu$ and $\nu^\prime$ are Poisson's ratios characterizing the lateral strain response in the plane of transverse isotropy to a stress acting parallel or normal to it, respectively, $G^\prime$ is the shear modulus in planes normal to the plane of transverse isotropy. Note also that $A_{11}-A_{12}=E/(1+\nu)$.

In the case of a transversely isotropic elastic layer bonded to a transversely isotropic elastic half-space, in accordance with the known solution \cite{Fabrikant2006} (see Eqs.~(102) and (103)), we will have
$$
\theta=\frac{1}{2\pi H^{(1)}},
$$
$$
L(u)=1+Q\bigl(\exp(-u/\gamma_1^{(1)}),\exp(-u/\gamma_2^{(1)})\bigr),
$$
where
$$
Q(x_1,y_1)=\frac{2M(x_1,y_1)}{N(x_1,y_1)}.
$$
According to Eq.~(101) \cite{Fabrikant2006}, the functions $M(x_1,y_1)$ and $N(x_1,y_1)$ are given by
\begin{eqnarray}
M(x_1,y_1) & = & -\bigl(g^{(1)}_1 a_{12}-g^{(1)}_2 a_{21}\bigr)x_1 y_1
-g^{(1)}_1 a_{11} x_1^2+g^{(1)}_2 a_{22} y_1^2
\nonumber \\
{} & {} & {}-(a_{12}a_{21}-a_{11} a_{22})x_1^2 y_1^2,
\nonumber \\
N(x_1,y_1) & = & 1+2\bigl(g^{(1)}_1 a_{12}-g^{(1)}_2 a_{21}\bigr)x_1 y_1
+\bigl(g^{(1)}_1+g^{(1)}_2\bigr)(a_{11} x_1^2-a_{22} y_1^2)
\nonumber \\
{} & {} & {}
+(a_{12}a_{21}-a_{11} a_{22})x_1^2 y_1^2.
\nonumber
\end{eqnarray}
According to Eqs.~(51)\,--\,(54) \cite{Fabrikant2006}, the coefficients $a_{11}$, $a_{12}$, $a_{21}$, and $a_{22}$ are defined as
\begin{eqnarray}
a_{11} & = & 1+\frac{2\gamma_1^{(1)}}{Z}\biggl(\biggl\{
-H^2\Bigl[\gamma_1^{(2)}-\bigl(m_1^{(2)}\bigr)^2\gamma_2^{(2)}\Bigr]
+gH\Bigl[2\bigl(\gamma_1^{(2)}-m_1^{(2)}\gamma_2^{(2)}\bigr)
\nonumber \\
{} & {} & {}+\bigl(m_1^{(1)}-1\bigr)\bigl(m_1^{(2)}-1\bigr)\gamma_2^{(1)}
\Bigr]\biggr\}\frac{1}{\gamma_1^{(2)}-\gamma_2^{(2)}}-g^2\biggr),
\nonumber \\
a_{12} & = & \frac{2\gamma_2^{(1)}}{Z}\biggl\{g^2 m_1^{(1)}
-gH\Bigl(g_1^{(2)}-m_1^{(2)}g_2^{(2)}\Bigr)\bigl(m_1^{(1)}+1\bigr)
\nonumber \\
{} & {} & {}+H^2\Bigl[g_1^{(2)}-\bigl(m_1^{(2)}\bigr)^2 g_2^{(2)}\Bigr]\biggr\},
\nonumber \\
a_{21} & = & \frac{2\gamma_1^{(1)}}{Z}\biggl\{g^2 m_1^{(1)}
-gH\Bigl(g_1^{(2)}-m_1^{(2)}g_2^{(2)}\Bigr)\bigl(m_1^{(1)}+1\bigr)
\nonumber \\
{} & {} & {}+H^2\Bigl[g_1^{(2)}-\bigl(m_1^{(2)}\bigr)^2 g_2^{(2)}\Bigr]\biggr\},
\nonumber \\
a_{22} & = & 1+\frac{2\gamma_2^{(1)}}{Z}\biggl(\biggl\{
H^2\Bigl[\gamma_1^{(2)}-\bigl(m_1^{(2)}\bigr)^2\gamma_2^{(2)}\Bigr]
-gH\Bigl[2m_1^{(1)}\bigl(\gamma_1^{(2)}-m_1^{(2)}\gamma_2^{(2)}\bigr)
\nonumber \\
{} & {} & {}-\bigl(m_1^{(1)}-1\bigr)\bigl(m_1^{(2)}-1\bigr)\gamma_1^{(1)}
\Bigr]\biggr\}\frac{1}{\gamma_1^{(2)}-\gamma_2^{(2)}}+g^2\bigl(m_1^{(1)}\bigr)^2
\biggr).
\nonumber
\end{eqnarray}
Here the following notation is used (see Eqs.~(48), (59), and (60) \cite{Fabrikant2006}):
\begin{eqnarray}
Z & = & H^2\Bigl\{
\bigl(\gamma_1^{(1)}-\gamma_2^{(1)}\bigr)
\Bigl(g_1^{(2)}-\bigl(m_1^{(2)}\bigr)^2 g_2^{(2)}\Bigr)\Bigr\}
\nonumber \\
{} & {} & {}-\frac{gH}{\gamma_1^{(2)}-\gamma_2^{(2)}}\Bigl[
\bigl(m_1^{(1)}-1\bigr)\bigl(m_1^{(2)}-1\bigr)
\Bigl(\gamma_1^{(1)}\gamma_2^{(1)}+\gamma_1^{(2)}\gamma_2^{(2)}\Bigr)
\nonumber \\
{} & {} & {}+2\Bigl(\gamma_1^{(2)}-m_1^{(2)}\gamma_2^{(2)}\Bigr)
\Bigl(\gamma_1^{(1)}-m_1^{(1)}\gamma_2^{(1)}\Bigr)\Bigr]
+g^2\Bigl[\gamma_1^{(1)}-\bigl(m_1^{(1)}\bigr)^2\gamma_2^{(1)}\Bigr],
\nonumber 
\end{eqnarray}
$$
g_1^{(1)}=\frac{\gamma_1^{(1)}}{\gamma_1^{(1)}-\gamma_2^{(1)}},\quad
g_2^{(1)}=\frac{\gamma_2^{(1)}}{\gamma_1^{(1)}-\gamma_2^{(1)}},\quad
g_1^{(2)}=\frac{\gamma_1^{(2)}}{\gamma_1^{(2)}-\gamma_2^{(2)}},\quad
g_2^{(2)}=\frac{\gamma_2^{(2)}}{\gamma_1^{(2)}-\gamma_2^{(2)}},
$$
$$
H=\frac{H^{(2)}\bigl(m_1^{(1)}-1\bigr)}{H^{(1)}\bigl(m_1^{(2)}-1\bigr)},\quad
g=\frac{\gamma_1^{(1)}-\gamma_2^{(1)}}{\gamma_1^{(2)}-\gamma_2^{(2)}}.
$$
We note two misprints in formula (59) \cite{Fabrikant2006} for $Z$: one closing bracket is missing and the product $gH$ is misprinted as $H_g$.

The layer material is characterized by the elastic constants $H^{(1)}$, $m_1^{(1)}$, $\gamma_1^{(1)}$, $\gamma_2^{(1)}$, and so on, while the material of the elastic half-space is characterized by the elastic constants $H^{(2)}$, $m_1^{(2)}$, $\gamma_1^{(2)}$, $\gamma_2^{(2)}$, and so on. So, we have
$$
H^{(n)}=\frac{\bigl(\gamma_1^{(n)}+\gamma_2^{(n)}\bigr)A_{11}^{(n)}}{
2\pi\Bigl(A_{11}^{(n)}A_{33}^{(n)}-\bigl(A_{13}^{(n)}\bigr)^2\Bigr)}
$$
and
$$
m_1^{(n)}=\frac{A_{11}^{(n)}\bigl(\gamma_1^{(n)}\bigr)^2-A_{44}^{(n)}}{A_{13}^{(n)}+A_{44}^{(n)}},\quad
m_2^{(n)}=\frac{A_{11}^{(n)}\bigl(\gamma_2^{(n)}\bigr)^2-A_{44}^{(n)}}{A_{13}^{(n)}+A_{44}^{(n)}}.
$$

Finally, the dimensionless parameters $\gamma_1^{(n)}$, $\gamma_2^{(n)}$ are defined as the roots of the equation 
\begin{equation}
\gamma^4 A_{11}^{(n)}A_{44}^{(n)}-\gamma^2\Bigl[A_{11}^{(n)}A_{33}^{(n)}-A_{13}^{(n)}
\Bigl(A_{13}^{(n)}+2A_{44}^{(n)}\Bigr)\Bigr]+A_{33}^{(n)}A_{44}^{(n)}=0.
\label{1ds(1.9)}
\end{equation}

In the case of isotropy, we have
$$
L(u)=1+Q(u),\quad
Q(u)=-2e^{-2u}\frac{d_1+d_2 e^{-2u}}{1+d_3 e^{-2u}+d_2 e^{-4u}},
$$
$$
\theta=\frac{E^{(1)}}{2\bigl(1-\bigl(\nu^{(1)}\bigr)^2\bigr)},
$$
where $E^{(1)}$ and $\nu^{(1)}$ are Young's modulus and Poisson's ratio of the layer material, while the coefficients $d_1$, $d_2$, and $d_3$ are given by formulas (105)\,--\,(108) \cite{Fabrikant2006}. We note that in this particular case, Fabrikant's solution obtained in \cite{Fabrikant2006} is in agreement with Burmister's solution \cite{Burmister1945a}.

\end{document}